\documentclass[a4paper]{amsart} 
\usepackage{amssymb} 
\usepackage{tikz-cd}
\usepackage[all]{xy}
\SelectTips{cm}{}
\usepackage{hyperref}
\calclayout

\title{Morphisms determined by objects and flat covers}
\author{Henning Krause}
\address{Henning Krause\\ Fakult\"at f\"ur Mathematik\\
  Universit\"at Bielefeld\\ D-33501 Bielefeld\\ Germany.}
\email{hkrause@math.uni-bielefeld.de}

\thanks{Version from June 25, 2014.}
\newtheorem{lemma}{Lemma}[section]
\newtheorem{proposition}[lemma]{Proposition}
\newtheorem{corollary}[lemma]{Corollary}
\newtheorem{theorem}[lemma]{Theorem}
\newtheorem{conjecture}[lemma]{Conjecture}

\theoremstyle{remark}
\newtheorem{remark}[lemma]{Remark}

\theoremstyle{definition}
\newtheorem{example}[lemma]{Example}

\newtheorem{problem}[lemma]{Problem}

\numberwithin{equation}{section}

\renewcommand{\mod}{\operatorname{\mathsf{mod}}\nolimits}
\DeclareMathOperator*{\colim}{colim}

\newcommand{\sub}{\operatorname{sub}\nolimits}

\newcommand{\rad}{\operatorname{rad}\nolimits}

\newcommand{\coind}{\operatorname{coind}\nolimits}
\newcommand{\ev}{\operatorname{ev}\nolimits}

\newcommand{\id}{\operatorname{id}\nolimits}
\newcommand{\Gr}{\operatorname{Gr}\nolimits}

\newcommand{\Mod}{\operatorname{\mathsf{Mod}}\nolimits}

\newcommand{\End}{\operatorname{End}\nolimits}
\newcommand{\fp}{\operatorname{\mathsf{fp}}\nolimits}
\newcommand{\Fp}{\operatorname{\mathsf{Fp}}\nolimits}
\newcommand{\Hom}{\operatorname{Hom}\nolimits}

\newcommand{\Image}{\operatorname{Im}\nolimits}
\newcommand{\Ker}{\operatorname{Ker}\nolimits}

\newcommand{\Ext}{\operatorname{Ext}\nolimits}

\newcommand{\Ab}{\mathsf{Ab}}
\newcommand{\op}{\mathrm{op}}

\newcommand{\fg}{\mathrm{fg}}

\newcommand{\lto}{\longrightarrow}
\newcommand{\xto}{\xrightarrow}

\def\a{\alpha}
\def\b{\beta}
\def\e{\varepsilon}

\def\p{\phi}

\def\Ga{\Gamma}
\def\La{\Lambda}

\def\A{{\mathsf A}}

\def\C{{\mathsf C}}

\def\Sc{{\mathsf S}}

\def\bbF{{\mathbb F}}

\def\bbQ{{\mathbb Q}}

\def\bbX{{\mathbb X}}
\def\bbZ{{\mathbb Z}}

\begin{document}

\begin{abstract}
  We describe a procedure for constructing morphisms in additive
  categories, combining Auslander's concept of a morphism determined
  by an object with the existence of flat covers. Also, we show how
  flat covers are turned into projective covers, and we interprete
  these constructions in terms of adjoint functors.
\end{abstract}

\maketitle

\section*{Introduction}

Functors and morphisms determined by objects were introduced by
Auslander in his Philadelphia notes \cite{Au1978}. These concepts
provide a method to construct and organise morphisms in additive
categories, generalising previous work of Auslander and Reiten on
almost split sequences \cite{AuRe1975}. More recently, Ringel
presented a survey of these results, rearranged them as lattice
isomorphisms (the Auslander bijections), and added a host of
interesting examples \cite{Ri2013}.

The starting point for our work is the following natural question: Is
there a procedure for constructing morphisms ending at a fixed object
in an additive category? More precisely, we are looking for
\begin{itemize}
\item[--]\emph{invariants} of morphisms ending at some fixed object,
  and 
\item[--] \emph{constructions} for universal morphisms with respect to
  these invariants.
\end{itemize}

An answer to this question is presented in Theorem~\ref{th:mor}. This
combines Auslander's concept of a morphism determined by an object
with a deep result about functor categories, which says that every
additive functor admits a flat cover \cite{BEE2001}.

The second part of this note is inspired by the first. We show in
Theorem~\ref{th:fc} that every flat cover is in fact a projective
cover, when viewed in an appropriate abelian category. Building on
another of Auslander's paradigms \cite{Au1966}, we use the Yoneda
embedding
\[\A\lto\Fp(\A^\op,\Ab),\quad X\mapsto\Hom_\A(-,X)\] 
of an additive category $\A$ into the category of finitely presented
functors $\A^\op\to\Ab$.  Take for instance the category $\A=\Mod \La$
of modules over a ring $\La$. It is somewhat surprising that any
functor $\A^\op\to\Ab$ preserving filtered colimits in $\A$ belongs to
$\Fp(\A^\op,\Ab)$ and admits a projective cover, even a minimal
projective presentation. This is precisely what we exploit in
Theorem~\ref{th:mor} and explain in Theorem~\ref{th:fc}.

The final section discusses the connection with almost split
sequences, and we complement Theorem~\ref{th:mor} by a non-existence
result for an almost split sequence ending at a module which is not
finitely presented (Proposition~\ref{pr:ass}).

\section{Morphisms determined by objects} 

\subsection*{Invariants of morphisms}
We fix an additive category $\A$ and let $\C$ be a set of objects. We
shall view $\C$ as a full subcategory of $\A$. A \emph{$\C$-module} is by
definition an additive functor $\C^\op\to\Ab$ into the category $\Ab$
of abelian groups, and a morphism between two $\C$-modules is a
natural transformation. The $\C$-modules form an abelian category
which we denote by $(\C^\op,\Ab)$. For example, if $\C$ consists of one
object $C$, then $(\C^\op,\Ab)$ is the category of modules over the
endomorphism ring of $C$. Note that (co)kernels and (co)products in
$(\C^\op,\Ab)$ are computed pointwise: for instance, a sequence $X\to Y\to
Z$ of morphisms between $\C$-modules is exact if and only if the
sequence $X(C)\to Y(C)\to Z(C)$ is exact in $\Ab$ for all $C$ in $\C$.

Every object $X$ in $\A$ gives rise to a $\C$-module \[\Hom_\A(\C, X)
= \Hom_\A(-, X)|_\C\colon\C^\op\lto\Ab\]
and every morphism $\a\colon X\to Y$ in $\A$ yields a morphism
\[\Hom_\A(\C,\a)\colon\Hom_\A(\C,X)\lto\Hom_\A(\C,Y).\] The image
$\Image\Hom_\A(\C,\a)$ of $\Hom_\A(\C,\a)$ is the invariant we shall
use. In fact, we construct a `right adjoint' which takes a submodule
of $\Hom_\A(\C,Y)$ to a morphism ending at $Y$.

\subsection*{Constructing morphisms}

We work in an additive category $\A$ which is \emph{locally finitely
  presented} \cite{CB1994}. This means the full subcategory $\fp\A$ of
finitely presented objects is essentially small and each object in
$\A$ is a filtered colimit of objects in $\fp\A$. Recall that an
object $X$ is \emph{finitely presented} if the functor
$\Hom_\A(X,-)\colon\A\to\Ab$ preserves filtered colimits. 

For example, the category $\Mod\La$ of modules over a ring $\La$ is
locally finitely presented. Then $\mod\La$ denotes the full
subcategory of finitely presented $\La$-modules.

A morphism $\a\colon X\to Y$ is called \emph{right minimal} if
every endomorphism $\p\colon X\to X$ with $\a\p=\a$ is invertible.

The following theorem is the main result of this work.

\begin{theorem}\label{th:mor} 
Let $\A$ be a locally finitely presented
additive category and $\C$ a set of finitely presented objects. For an
object $Y$ in $\A$ and a submodule $H\subseteq \Hom_\A(\C,Y)$, there
exists a morphism $\a\colon X\to Y$ in $\A$ (unique up to non-unique
isomorphism) such that the following holds:
\begin{enumerate}
\item $\Image\Hom_\A(\C,\a)=H$ and any morphism $\a'\colon X'\to Y$
with $\Image\Hom_\A(\C,\a')\subseteq H$ factors through $\a$.
\item $\a$ is right minimal.
\end{enumerate}
\end{theorem}

The proof will be given later in this section. A second and more
elementary proof for $\A=\Mod\La$ can be found in an appendix.

\subsection*{Morphisms determined by objects} 

Following \cite{Au1978}, a morphism $\a\colon X\to Y$ in $\A$ is
called \emph{right determined by $\C$} (or simply \emph{right
  $\C$-determined}) if any morphism $\a'\colon X'\to Y$ satisfying
\[\Image\Hom_\A(\C,\a')\subseteq \Image\Hom_\A(\C,\a)\] factors through
$\a$. In fact, Auslander established Theorem~\ref{th:mor} for module
categories with $\C$ consisting of a single object
\cite[Theorem~I.3.19]{Au1978}, generalising previous work of Auslander
and Reiten on almost split sequences \cite{AuRe1975}.

An obvious question to ask is when a morphism is right determined by
some set of finitely presented objects. We give a general answer in
Proposition~\ref{pr:kernel2}.  For an Artin algebra, every morphism
between finitely presented modules is right determined by some
finitely presented module \cite[Theorem~2.6]{Au1978b}.

\subsection*{Functoriality}

The assignment $(\C, H,Y)\mapsto \a$ in Theorem~\ref{th:mor} is
functorial. This has been pointed out by Ringel in his survey of
Auslander's results \cite{Ri2013}.

Let us formulate the functoriality. We fix an object $Y\in\A$.  The
morphisms ending at $Y$ are \emph{preordered}. This means for
morphisms $\a\colon X\to Y$ and $\a'\colon X'\to Y$ that $\a'\leq\a$
when $\a'$ factors through $\a$. We obtain a \emph{poset} by
identifying $\a$ and $\a'$ when $\a'\leq\a$ and $\a\leq\a'$. Let us
denote this poset by $[\A/Y]$ because it is derived from the slice
category $\A/Y$.\footnote{In \cite{Ri2013}, the poset $[\A/Y]$ is
  called \emph{right factorisation lattice} for $Y$ and is denoted by
  $[\to Y\rangle$.}

For a pair $\C\subseteq\fp\A$ and $H\subseteq\Hom_\A(\C,Y)$ let
$\a_{\C,H}\colon X_{\C,H}\to Y$ denote the right minimal and right
$\C$-determined morphism such that $\Image\Hom_\A(\C,\a_{\C,H})=H$; it
exists and is well-defined up to a non-unique isomorphism by
Theorem~\ref{th:mor}. Note that $\a_{\C,H}$ is unique when viewed as
an element of $[\A/Y]$.

\begin{lemma}
  Let $\phi\colon Y'\to Y$ be a morphism in $\A$, $\C\subseteq \C'\subseteq\fp\A$,
  $H\subseteq\Hom_\A(\C,Y)$, and $H'\subseteq\Hom_\A(\C',Y')$. Then
  $\Hom_\A(\C,\phi)(H')\subseteq H$ implies  $\phi\a_{\C',H'}\leq\a_{\C,H}$.
\end{lemma}
\begin{proof}
The assumptions imply  
\[\Image\Hom_\A(\C,\phi\a_{\C',H'})=
\Hom_\A(\C,\phi)(\Image\Hom_\A(\C,\a_{\C',H'})) \subseteq H.\] Thus
$\phi\a_{\C',H'}$ factors through $\a_{\C,H}$.
\end{proof}

For an object $X$ in an abelian category let $\sub(X)$ denote its
lattice of subobjects.

\begin{remark}
Viewing a poset as a category, the map
\[\sub(\Hom_\A(\C,Y))\lto [\A/Y],\quad H\mapsto \a_{\C,H}\]  
is right adjoint to
\[ [\A/Y] \lto \sub(\Hom_\A(\C,Y)),\quad \a\mapsto \Image\Hom_\A(\C,\a).\] 
\end{remark}

There are some natural choices of triples $(\C,H,Y)$. For example,
right almost split morphisms arise from triples $(Y,\rad\End_\A(Y),Y)$
for a finitely presented object $Y$ with local endomorphism ring
\cite[\S II.2]{Au1978}. Another obvious choice is $H=0$; not much
seems to be known in this case.

\begin{problem}[Auslander \cite{Au1992}]
Describe the morphism $\a_{\C,0}$  for $Y\in\A$ and $\C\subseteq\fp\A$.
\end{problem}

Note the extremes: $\a_{\varnothing,0}=\id_Y$ is the identity morphism and
$\a_{\C,0}=0$ when $\C$ contains a generator of $\fp\A$.

The assignment $(\C,H)\mapsto\a_{\C,H}$ has been studied in some
detail for modules over Artin algebras. We include Ringel's
formulation of the Auslander bijection as an example.

\begin{example}[Ringel \cite{Ri2013}]
  Let $\La$ be an Artin algebra. For $Y\in\mod\La$ the assignment
  $(\C,H)\mapsto\a_{\C,H}$ induces an isomorphism\footnote{The colimit
    is taken over the collection of finite subsets
    $\C\subseteq\mod\La$, identifying $\C=\{C_1,\ldots,C_n\}$ with
    $C=C_1\oplus\ldots\oplus C_n$.}
  \[\colim_{C\in\mod\La}\sub(\Hom_\La(C,Y))\stackrel{\sim}\lto [\mod\La/Y].\]
\end{example}

\subsection*{Functors determined by objects}

Let $\A$ be an additive category. We consider additive functors
$\A^\op\to\Ab$ into the category of abelian groups; they form an
abelian category which we denote by $(\A^\op,\Ab)$. 

Fix an additive functor $F\colon\A^\op\to\Ab$ and a set $\C$ of
objects in $\A$ (viewed as a full subcategory). We write $F|_\C$ for
the restriction $\C^\op\to\Ab$.  Following Auslander \cite{Au1978}, a
subfunctor $F'\subseteq F$ is called \emph{$\C$-determined} if for any
subfunctor $F''\subseteq F$
  \[F''\subseteq F'\quad\iff\quad F''|_\C\subseteq F'|_{\C}.\] 

For a subfunctor  $H\subseteq F|_\C$ define  a subfunctor
$F_H\subseteq F$ by
\begin{equation}\label{eq:det}
F_H(X)=\bigcap_{\genfrac{}{}{0pt}{}{\a\colon C\to X}{C\in\C}} F(\a)^{-1}(H(C))\qquad
\text{for }X\in\A.
\end{equation}

\begin{lemma}[{Auslander \cite[Proposition~I.1.2]{Au1978}}]\label{le:det}
The subfunctor
 $F_H\subseteq F$ is $\C$-determined and $F_H|_\C=H$. 
\end{lemma}
\begin{proof}
Let $F'\subseteq F$ be a subfunctor and $F'|_\C\subseteq H$. This
implies for each morphism $\a\colon C\to X$ with $C\in\C$ that
$F'(X)\subseteq  F(\a)^{-1}(H(C))$. Thus $F'\subseteq F_H$.
\end{proof}

\begin{remark}
  Let $F\in(\A^\op,\Ab)$. Viewing a poset as a category, the map
\[\sub (F|_\C)\lto\sub(F),\quad H\mapsto F_H\] is right adjoint to
\[\sub (F)\lto\sub(F|_\C),\quad G\mapsto G|_\C.\]
\end{remark}

We include an explicit example.

\begin{example}
Let $\La$ be a ring and  consider the forgetful functor
$F\colon\mod\La\to \Ab$ which takes a module to its underlying abelian group.

(1)  For $C\in\mod\La$ and an $\End_\La(C)$-submodule $H\subseteq C$, the
$C$-determined subfunctor $F_H$ is given by
\[F_H(X)=\bigcap_{\a\colon X\to C}\{ x\in X\mid \a(x)\in
H\}\qquad\text{for }X\in\mod\La.\] 

(2) For an Artin algebra $\La$, let $\sub_\fg(F)$ denote the poset of
finitely generated subfunctors of $F$. Then the assignment $H\mapsto
F_H$ induces an isomorphism
\[\colim_{C\in \mod\La}\sub(C)\stackrel{\sim}\lto\sub_\fg(F).\]

(3) Consider the algebra $\La=\bbF_2[\e]$ of dual numbers ($\e^2=0$)
and the $\La$-module $C=\bbF_2[\e]\oplus \bbF_2$. Then the poset
$\sub(C)$ is isomorphic to $\sub_\fg(F)$ and its Hasse diagram is the
following.
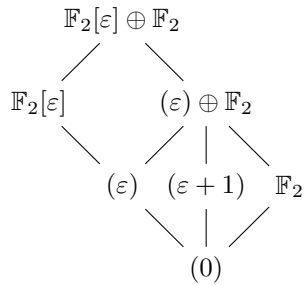
\[\begin{tikzpicture}[scale=.55]
  \node (a) at (-1,3) {$\bbF_2[\e]\oplus \bbF_2$};
  \node (b) at (-3,1) {$\bbF_2[\e]$};
  \node (c) at (1,1) {$(\e)\oplus \bbF_2$};
  \node (d) at (-1,-1) {$(\e)$};
  \node (e) at (1,-1) {$(\e+1)$};
  \node (f) at (3,-1) {$\bbF_2$};
  \node (g) at (1,-3) {$(0)$};
  \draw (d) -- (c) -- (a) -- (b) -- (d) -- (g) -- (e) -- (c) -- (f) -- (g);
\end{tikzpicture}
\]
\end{example}

\subsection*{Flat functors and flat covers}
Let $\C$ be an essentially small additive category. We consider the
category $(\C^\op,\Ab)$ of additive functors
$F\colon\C^\op\to\Ab$. Recall that $F$ is \emph{flat} if it can be
written as a filtered colimit of representable functors.

The following result establishes a connection between locally
finitely presented additive categories and categories of flat functors.

\begin{theorem}[{Crawley-Boevey \cite[\S1.4]{CB1994}}]\label{th:flat}
Let $\A$ be a  locally finitely presented category. Then the Yoneda functor
  \[h\colon\A\lto ((\fp\A)^\op,\Ab),\quad X\mapsto
  \Hom_\A(-,X)|_{\fp\A}\] identifies $\A$ with the full subcategory of
  flat functors $(\fp\A)^\op\to\Ab$.\qed 
\end{theorem}

A morphism $\pi\colon F\to G$ in $(\C^\op,\Ab)$ is a \emph{flat
  cover} of $G$ if the following holds:
\begin{enumerate}
\item $F$ is flat and every morphism $F'\to G$ with $F'$ flat factors
  through $\pi$.
\item $\pi$ is right minimal.
\end{enumerate}
A \emph{minimal flat presentation} of $G$ is an exact sequence
\[F_1\lto F_0\stackrel{\pi}\lto G\lto 0\] such that $F_0\to G$ and
$F_1\to\Ker\pi$ are flat covers.  A \emph{projective cover} and a
\emph{minimal projective presentation} are defined analogously,
replacing the term flat by projective.\footnote{This definition of a
  projective cover is equivalent to the usual one
  which requires the kernel to be superfluous.}

\begin{theorem}[{Bican--El Bashir--Enochs \cite{BEE2001}}]\label{th:flatcover}
Every additive functor $\C^\op\to\Ab$ admits a flat cover.\qed
\end{theorem}

\subsection*{Proof of the main theorem}

We are ready to prove  Theorem~\ref{th:mor}. The basic idea is to
identify $\A$ with the category of flat functors $(\fp\A)^\op\to\Ab$
and to employ the existence of flat covers.

\begin{proof}[Proof of Theorem~\ref{th:mor}]
  We apply Theorems~\ref{th:flat} and \ref{th:flatcover}.  Consider
  the subfunctor $H\subseteq h(Y)|_\C$ and choose a flat cover
  $\pi\colon h(X)\to h(Y)_H$ in $((\fp\A)^\op,\Ab)$. The composite
  $h(X)\to h(Y)_H\to h(Y)$ is of the form $h(\a)$ for some morphism
  $\a\colon X\to Y$ in $\A$. We check the properties of $\a$ and apply
  Lemma~\ref{le:det}. Thus $\Image\Hom_\A(\C,\a)=h(Y)_H|_\C=H$. For a
  morphism $\a'\colon X'\to Y$ in $\A$ with
  $\Image\Hom_\A(\C,\a')\subseteq H$, we have $\Image
  h(\a')|_\C\subseteq H$ and therefore $\Image h(\a')\subseteq h(Y)_H=
  \Image h(\a)$. Using that $\pi$ is a flat cover, we obtain a
  morphism $\p\colon X'\to X$ in $\A$ satisfying $\a\p=\a'$. The fact
  that $\pi$ is right minimal implies that $\a$ is right minimal.
\end{proof}

\begin{corollary}\label{co:sum}
  Let $\a\colon X\to Y$ be a morphism which is right determined by a
  set of finitely presented objects. Then there is an essentially
  unique decomposition $X=X'\oplus X''$ such that $\a|_{X'}$ is right
  minimal and $\a|_{X''}=0$.
\end{corollary}
\begin{proof}
  Suppose that $\a$ is right $\C$-determined and let
  $H=\Image\Hom_\A(\C,\a)$.  Then there exists a right minimal and
  right $\C$-determined morphism $\a'\colon X'\to Y$ with
  $\Image\Hom_\A(\C,\a') =H$ by Theorem~\ref{th:mor}. The minimality
  of $\a'$ yields a decomposition $X=X'\oplus X''$ such that
  $\a|_{X'}=\a'$ and $\a|_{X''}=0$.
\end{proof}

\subsection*{Weak kernels}

Let $\A$ be an additive category. A morphism $X\to Y$ is a \emph{weak
  kernel} of a morphism $Y\to Z$ in $\A$ if the induced
sequence \[\Hom_\A(-,X)\lto \Hom_\A(-,Y)\lto \Hom_\A(-,Z)\] is
exact. A weak kernel is called \emph{minimal} if it is a right minimal
morphism. Note that a minimal weak kernel is unique up to a non-unique
isomorphism; it is a kernel if a kernel exists.

\begin{proposition}\label{pr:kernel}
In a locally finitely presented additive category every morphism
admits a minimal weak kernel.
\end{proposition}
\begin{proof}
  We apply Theorems~\ref{th:flat} and \ref{th:flatcover}.  Fix a
  morphism $\b\colon Y\to Z$ and choose a flat cover $h(X)\to \Ker
  h(\b)$ in $((\fp\A)^\op,\Ab)$. The composite $h(X)\to \Ker h(\b)\to
  h(Y)$ is of the form $h(\a)$ for some morphism $\a\colon X\to Y$ in
  $\A$, which is a minimal weak kernel for $\b$.
\end{proof}

\subsection*{Morphisms determined by finitely presented objects} 

Let $\A$ be a locally finitely presented additive category.  Recall
that a sequence $0\to X\to Y\to Z\to 0$ of morphisms in $\A$ is
\emph{pure-exact} if for each $C\in\fp\A$ the induced sequence
\[0\to\Hom_\A(C,X)\to\Hom_\A(C,Y)\to\Hom_\A(C,Z)\to 0\] is exact. An
object $X$ in $\A$ is \emph{pure-injective} if every pure-exact
sequence $0\to X\to Y\to Z\to 0$ is split exact.

The following proposition characterises the morphisms which are right
determined by a set of finitely presented objects.  The relevance of
pure-injectives was noticed by Auslander in some special cases
\cite[\S I.10]{Au1978}.

\begin{proposition}\label{pr:kernel2}
  Let $\A$ be a locally finitely presented additive category. 
\begin{enumerate}
\item If a morphism in $\A$ is right minimal and right determined by a
  set of finitely presented objects in $\A$, then its minimal weak
  kernel is pure-injective.
\item If a morphism in $\A$ has a pure-injective kernel, then it is
  right determined by a set of finitely presented objects in $\A$.
\end{enumerate}
\end{proposition}
\begin{proof}
  We set $\C=\fp\A$ and recall that a functor $F\in (\C^\op,\Ab)$ is
  \emph{cotorsion} if $\Ext^1(E,F)=0$ for each flat $E\in
  (\C^\op,\Ab)$. Clearly, $X\in\A$ is pure-injective if and only if
  $h(X)$ is cotorsion.

  Observe that the kernel of a flat cover in $(\C^\op,\Ab)$ is
  cotorsion. Also, the flat cover of a cotorsion functor is cotorsion;
  see \cite[\S2]{BEE2001} for details.

  Now fix a sequence $X\xto{\a} Y\xto{\b} Z$ of morphisms such that
  $\a$ is the minimal weak kernel of $\b$. If $\b$ is right minimal
  and right determined by any set of finitely presented objects, then
  $h(Y)\to \Image h(\b)$ is a flat cover; this follows from the proof
  of Theorem~\ref{th:mor}.  Thus $\Image h(\a)=\Ker h(\b)$ is
  cotorsion. The proof of Proposition~\ref{pr:kernel} shows that
  $h(X)\to \Image h(\a)$ is a flat cover. Thus $h(X)$ is cotorsion and
  therefore $X$ is pure-injective. Conversely, if $\a$ is a kernel of
  $\b$ and $X$ is pure-injective, then every morphism $h(Y')\to \Image
  h(\b)$ factors through $h(Y)\to \Image h(\b)$. Thus every morphism
  $\b'\colon Y'\to Z$ with
  $\Image\Hom_\A(\C,\b')\subseteq\Image\Hom_\A(\C,\b)$ factors through
  $\b$.  This means that $\b$ is right $\C$-determined.
\end{proof}

In view of Corollary~\ref{co:sum} we obtain the following consequence.

\begin{corollary}
  Let $\a\colon X\to Y$ be a morphism which admits a kernel. Then $\a$
  is right determined by a set of finitely presented objects if and
  only if there is a decomposition $X=X'\oplus X''$ such that
  $\Ker\a|_{X'}$ is pure-injective and $\a|_{X''}=0$.\qed
\end{corollary}

\section{Flat versus projective covers}

In this section we show how flat covers are turned into projective
covers.  In fact, we prove a general result about locally finitely
presented additive categories which is tantamount to the the existence
of flat covers in functor categories.

\subsection*{Finitely presented functors}

Let $\A$ be an additive category. We denote by $\Fp(\A^\op,\Ab)$ the
category of finitely presented functors $F\colon\A^\op\to\Ab$. Recall
that $F$ is \emph{finitely presented} (or \emph{coherent}) if it fits
into an exact sequence
\[\Hom_\A(-,X)\lto \Hom_\A(-,Y)\lto F\lto 0.\]
We call such a presentation \emph{minimal} if the morphisms
$\Hom_\A(-,Y)\to F$ and $\Hom_\A(-,X)\to \Hom_\A(-,Y)$ are right
minimal.

The following lemma is well-known and easily proved (using Yoneda's lemma).

\begin{lemma}\label{le:fp}
  The category $\Fp(\A^\op,\Ab)$ is abelian if and only $\A$ admits
  weak kernels. In this case the Yoneda embedding
  $X\mapsto\Hom_\A(-,X)$ identifies the idempotent completion of $\A$
  with the category of projective objects in $\Fp(\A^\op,\Ab)$.\qed
\end{lemma}

If $\A$ is locally finitely presented then we consider the \emph{evaluation functor}
\begin{equation}\label{eq:ev}
\Fp(\A^\op,\Ab)\lto ((\fp\A)^\op,\Ab),\quad F\mapsto F|_{\fp\A}.
\end{equation}
The following theorem establishes a right adjoint which identifies
flat covers with projective covers.

\begin{theorem}\label{th:fc}
Let $\A$ be a locally finitely presented category. 
\begin{enumerate}
\item The category $\Fp(\A^\op,\Ab)$ of finitely presented functors
  $\A^\op\to\Ab$ is abelian.
\item For an additive functor $F\colon (\fp\A)^\op\to\Ab$, the unique
  functor $\widetilde F\colon\A^\op\to\Ab$ extending $F$ and
  preserving filtered colimits in $\A$ is finitely presented and admits a minimal
  projective presentation in $\Fp(\A^\op,\Ab)$.
\item The assignment $F\mapsto\widetilde F$ provides a fully faithful
  right adjoint to the evaluation functor \eqref{eq:ev}.
\end{enumerate}
\end{theorem}

We postpone the proof and illustrate the theorem by a couple of
examples. The first one shows how flat covers of modules over a ring
are derived from this result.

\begin{example}
  Let $\La$ be a ring and denote by $\A$ the category of flat
  $\La$-modules. Then $\fp\A$ equals the category of finitely
  generated projective $\La$-modules and evaluation at $\La$ yields an
  equivalence $((\fp\A)^\op,\Ab)\xto{\sim}\Mod\La$ (which we view as
  identification). For a $\La$-module $Y$, the theorem yields a
  projective cover $\Hom_\A(-,X)\to\widetilde Y$. Evaluation at
  $\La$ then gives a flat cover $X\to Y$.
\end{example}

\begin{example}
  Not all functors in $\Fp(\A^\op,\Ab)$ admit a projective cover. For
  instance, take $\A=\Mod\bbZ$ and consider the canonical morphism
  $\bbZ\to\bbZ/p$ for any prime $p$. Then the image of the induced
  morphism $\Hom_\bbZ(-,\bbZ)\to\Hom_\bbZ(-,\bbZ/p)$ admits no
  projective cover, because a projective cover $\Hom_\bbZ(-,X)\to F$
  would give a projective cover $X\to\bbZ/p$ in $\A$ (which is known
  not to exist).
\end{example}

\begin{proof}[Proof of Theorem~\ref{th:fc}]
  The category $\A$ has weak kernels by
  Proposition~\ref{pr:kernel}. Thus  $\Fp(\A^\op,\Ab)$ is abelian by Lemma~\ref{le:fp}.

  Now fix $F\in ((\fp\A)^\op,\Ab)$ and observe that $\widetilde
  F\cong\Hom(h-,F)$ since $h$ preserves filtered colimits. This
  yields for $X\in\A$ a  functorial isomorphism
\begin{equation}\label{eq:adj}
\Hom(\Hom_\A(-,X),\widetilde F)\cong \widetilde F(X)\cong\Hom(\Hom_\A(-,X)|_{\fp\A},F)
\end{equation}
which extends to an isomorphism
\[\Hom(E,\widetilde F)\cong \Hom(E|_{\fp\A},F)\] for all
$E\in\Fp(\A^\op,\Ab)$.  The adjointness property of the assignment
$F\mapsto \widetilde F$ then follows. Also, the functor is fully
faithful since $\widetilde F|_{\fp\C}=F$, and it identifies the flat
functors in $((\fp\A)^\op,\Ab)$ with the projective objects in
$\Fp(\A^\op,\Ab)$.

Next we show that $\widetilde F$ is finitely presented. In fact, we
obtain a minimal projective presentation of $\widetilde F$ in
$\Fp(\A^\op,\Ab)$ by choosing a minimal flat presentation
\[h(X)\lto h(Y)\stackrel{\pi}\lto F\lto 0\] in $((\fp\A)^\op,\Ab)$;
see Theorem~\ref{th:flatcover}. From \eqref{eq:adj} it follows that
the corresponding sequence
\[\Hom_\A(-,X)\lto \Hom_\A(-,Y)\stackrel{\widetilde\pi}\lto \widetilde
F\lto 0\] is a minimal projective presentation.
\end{proof}

\begin{corollary}\label{co:fc}
  The evaluation functor \eqref{eq:ev} is exact and admits both
  adjoints; it induces an equivalence
\[\frac{\Fp(\A^\op,\Ab)}{\{F\mid F|_{\fp\A}=0 \}}\stackrel{\sim}\lto ((\fp\A)^\op,\Ab).\]
\end{corollary}
\begin{proof}
  The right adjoint exists by Theorem~\ref{th:fc}, and the left
  adjoint is the unique colimit preserving functor sending
  $\Hom_\A(-,X)$ to $\Hom_\A(-,X)$ for all $X\in\fp\A$. 

  Let $\Sc$ denote the kernel of the evaluation functor.  A
  quasi-inverse for the functor $\Fp(\A^\op,\Ab)/\Sc\to
  ((\fp\A)^\op,\Ab)$ is obtaind by composing the left (or right)
  adjoint of the evaluation functor with the quotient functor
  $\Fp(\A^\op,\Ab)\to \Fp(\A^\op,\Ab)/\Sc$.
\end{proof}

\subsection*{Auslander's formula}
Let $\A$ be an abelian category.  A somewhat hidden result in
Auslander's account on coherent functors \cite[p.~205]{Au1966} shows
that the Yoneda functor $\A\to\Fp(\A^\op,\Ab)$ admits an exact left
adjoint which sends a representable functor $\Hom_\A(-,X)$ to $X$ and
yields \emph{Auslander's formula} \cite{Le1998}
\[\frac{\Fp(\A^\op,\Ab)}{\{F\mid F \text{ is exactly presented} \}}\stackrel{\sim}\lto \A.\]
Here, a functor $F\colon\A^\op\to\Ab$ is
\emph{exactly presented} if it fits into an exact sequence
\begin{equation}\label{eq:fp}
0\to\Hom_\A(-,X)\to \Hom_\A(-,Y)\to \Hom_\A(-,Z)\to F\to 0
\end{equation} 
such that the corresponding sequence $0\to X\to Y\to Z\to 0$ in $\A$
is exact. 

One may think of Auslander's formula as a prototype for
Corollary~\ref{co:fc}. In particular, we see that the kernel of the
evaluation functor is given by the functors $F\colon\A^\op\to\Ab$ with
a presentation \eqref{eq:fp} such that the corresponding sequence
$0\to X\to Y\to Z\to 0$ in $\A$ is pure-exact.

\section{Almost split sequences}

The concept of a morphism determined by an object generalises that of
an almost split morphism.  In fact, a morphism $\a\colon X\to Y$ in an
additive category $\A$ is \emph{right almost split} (that is, $\a$ is
not a retraction and every morphism $X'\to Y$ that is not a retraction
factors through $\a$) if and only if $\Ga=\End_\A(Y)$ is a local ring,
$\a$ is right determined by $Y$, and $\Image\Hom_\A(Y,\a)=\rad\Ga$
\cite[\S II.2]{Au1978}.

Now let $\A$ be an abelian category. Recall that an exact sequence $0\to
X\xto{\a}Y\xto{\b} Y\to 0$ is \emph{almost split} if the morphism $\a$
is left almost split and $\b$ is right almost split.  This is
equivalent to $\b$ being a right minimal and right almost split
morphism \cite[\S II.4]{Au1978}.

The construction of morphisms in Theorem~\ref{th:mor} requires the
determining objects to be finitely presented. The following
proposition shows that this assumption is necessary. A similar
non-existence result for finite dimensional representations of
infinite quivers is due to Paquette \cite{Pa2012}.

\begin{proposition}\label{pr:ass}
The category of modules over $\bbZ$ does not admit an almost split
sequence $0\to X\to Y\to Z\to 0$ for $Z=\bbQ$.
\end{proposition}

Thus the category $\Mod\bbZ$ does not admit a right minimal and right
determined morphism for the triple $(\C,H,Y)=(\bbQ,0,\bbQ)$.

The end terms of an almost split sequence determine each other, and
this correpondence enjoys some weak functoriality. The proof of
Proposition~\ref{pr:ass} uses this argument, and I am grateful to
Helmut Lenzing for suggesting it.

\begin{lemma}[{\cite[Lemma~A.10]{Kr2012}}]\label{le:ass}
  \pushQED{\qed} Let $0\to X\to Y\to Z\to 0$ be an almost split
  sequence in any abelian category $\A$. Then there is an
  isomorphism \[\End_\A(X)/\rad\End_\A(X)
  \cong\End_\A(Z)/\rad\End_\A(Z).\qedhere\]
\end{lemma}
\begin{proof}[Proof of Proposition~\ref{pr:ass}]
  Suppose there exists an almost split sequence $0\to X\to Y\to
  \bbQ\to 0$ in $\Mod\bbZ$.  Then the lemma implies
\[\End_\bbZ(X)/\rad \End_\bbZ(X)\cong\bbQ.\] 
It follows that $X$ is divisible since multiplication with any
non-zero integer induces an endomorphism of $X$ which is invertible.
Thus the sequence splits which is impossible.
\end{proof}

Very little seems to be known about the non-existence of almost split
sequences. We conjecture the following.

\begin{conjecture}
Given an almost split sequence $0\to X\to Y\to Z\to 0$ in any module
category, the module $Z$ is finitely presented.
\end{conjecture}

\begin{appendix}

\section{Functors on module categories}

The proof of Theorem~\ref{th:mor} uses the existence of flat
covers. In this appendix we provide an elementary argument for the
case of a module category $\A=\Mod\La$. We proceed in three steps.

\subsection*{Minimal presentations}
The existence of flat covers in $((\mod\La)^\op,\Ab)$ can be
phrased as follows.

\begin{theorem}\label{th:cover}
  Let $\La$ be a ring. An additive functor
  $F\colon(\Mod\La)^\op\to\Ab$ preserving filtered colimits in
  $\Mod\La$ admits a minimal presentation
\[\Hom_\La(-,X)\lto \Hom_\La(-,Y)\lto F\lto 0.\]
\end{theorem}

\begin{remark}
(1) The evaluation functor $\Fp((\Mod\La)^\op,\Ab)\to ((\mod\La)^\op,\Ab)$
sends a minimal projective presentation to a minimal flat
presentation.

(2) The theorem complements a result of Crawley-Boevey \cite{CB1998}:
An additive functor $F\colon\Mod\La\to\Ab$ preserves filtered colimits
and products if and only if it admits a presentation
\[\Hom_\La(Y,-)\lto \Hom_\La(X,-)\lto F\lto 0\]
with $X,Y\in\mod\La$.
\end{remark}

\begin{proof}
  Any $\La$-module can be written as a filtered colimit of
  finitely presented modules. It follows that $F$ is determined by its
  restriction $F|_{\mod\La}$.  The functors \[I_C=
  \Hom_\bbZ(\Hom_\La(C,-),\bbQ/\bbZ)\] with $C\in\mod\La$ form a set
  of injective cogenerators for $((\mod\La)^\op,\Ab)$ and this yields
  an injective copresentation
\[0\lto F|_{\mod\La} \lto\prod_i I_{C_i}\lto\prod_j I_{C_j}\]
with $C_i,C_j\in\mod\La$. Thus we obtain an exact sequence
\[0\lto F \lto\prod_i I_{C_i}\lto\prod_j I_{C_j}\] in
$((\Mod\La)^\op,\Ab)$. It suffices to show for each $C\in\mod\La$
that the functor $I_C$ is finitely presented. Then one uses that the finitely presented
functors are closed under taking products and kernels.  Choose a
presentation
\[P_1\lto P_0\lto C\lto 0\] such that each $P_i$ is finitely generated
projective. This yields a presentation
\[I_{P_1}\lto I_{P_0}\lto I_{C}\lto 0.\]
We have
\begin{align*}
\Hom_\bbZ(\Hom_\Lambda(P_i,-),\bbQ/\bbZ) 
&\cong\Hom_\bbZ(-\otimes_\Lambda P_i^*,\bbQ/\bbZ) \\
&\cong\Hom_\La(-,\Hom_\bbZ(P_i^*,\bbQ/\bbZ))
\end{align*}
where $P_i^*=\Hom_\La(P_i,\La)$. Thus $I_C$ is finitely presented.

The morphisms $\Hom_\La(-,Y)\to F$ and $\Hom_\La(-,X)\to\Hom_\La(-,Y)$
can be chosen to be right minimal; this follows from a standard
argument \cite[\S7]{En1981}.
\end{proof}

\subsection*{Coinduction}

Let $\A$ be an additive category. For a full subcategory
$\C\subseteq\A$ consider the \emph{evaluation functor}
\[\ev_\C\colon (\A^\op,\Ab)\lto (\C^\op,\Ab),\quad F\mapsto F|_\C\]
and its right adjoint, the \emph{coinduction functor}
\[\coind_\C\colon (\C^\op,\Ab)\lto (\A^\op,\Ab)\]
given by \[\coind_\C I(X)=\Hom(\Hom_\A(-,X)|_\C,I)\qquad \text{for }I\in (\C^\op,\Ab),\,X\in\A.\]
For $F\in (\A^\op,\Ab)$ and $I\in (\C^\op,\Ab)$, the isomorphism
\begin{equation}\label{eq:Yon}
\Hom(F,\coind_\C I) \xto{\sim}\Hom(F|_\C,I),
\quad \eta\mapsto \bar\eta
\end{equation}
is given by $\bar\eta=\eta|_\C$, where we identify $(\coind_\C
I)|_\C\xto{\sim} I$.

The coinduction functor assists in understanding the assignment \eqref{eq:det}
\[F|_\C\supseteq H\longmapsto F_H\subseteq F\] for an additive functor
$F\colon\A^\op\to\Ab$.

\begin{lemma}\label{le:detfun}
  Let $F\colon\A^\op\to\Ab$ be an additive functor and $H\subseteq
  F|_\C$ a subfunctor.  For a morphism $\bar\eta\colon F|_\C\to I$
  with kernel $H$ we have $F_H=\Ker\eta$, where $\eta$ and $\bar\eta$
  are related via \eqref{eq:Yon}.
\end{lemma}
\begin{proof}
  Let $F'\subseteq F$ be a subfunctor. The isomorphism \eqref{eq:Yon}
  is functorial and the inclusion $F'\to F$ gives $\eta|_{F'}=0$ iff
  $\bar\eta|_{F'|_\C}=0$. Thus $F'\subseteq \Ker\eta$ iff
  $F'|_\C\subseteq H$. This means that $\Ker\eta$ is $\C$-determined,
  and $\Ker\eta|_\C=H$ by construction.
\end{proof}

\subsection*{Morphisms determined by objects}

The following lemma provides the connection between morphisms
and functors determined by objects.

\begin{lemma}\label{le:rightdet}
  Let $\A$ be an additive category and $\C\subseteq\A$.  A morphism
  $\a\colon X\to Y$ in $\A$ is right $\C$-determined if and only if
  the subfunctor
\[\Image\Hom_\A(-,\a)\subseteq\Hom_\A(-,Y)\] is $\C$-determined. 
\end{lemma}

\begin{proof} 
  This is clear because a morphism $\a'\colon X'\to Y$ factors through
  $\a$ iff\[\Image\Hom_\A(-,\a')\subseteq
  \Image\Hom_\A(-,\a).\qedhere\]
\end{proof} 

We are now ready for a second proof of Theorem~\ref{th:mor} for $\A=\Mod\La$.

\begin{proof}[Proof of Theorem~\ref{th:mor}]
  Set $F=\Hom_\A(-,Y)$. The subfunctor $F_H$ arises as a kernel of a
  morphism $F\to\coind_\C I$ for some $I\in (\C^\op,\Ab)$ by
  Lemma~\ref{le:detfun}. The functor $\coind_\C I$ preserves filtered
  colimits in $\A$ since $\C$ consists of finitely presented
  objects. Also $F$ preserves filtered colimits in $\A$, and therefore
  $F_H$ has this property. Thus $F_H$ admits a projective cover
  $\pi\colon\Hom_\A(-,X)\to F_H$ by Theorem~\ref{th:cover}. The
  composite $\Hom_\A(-,X)\to F_H\to F$ is represented by a morphism
  $\a\colon X\to Y$ which is right $\C$-determined by
  Lemma~\ref{le:rightdet}, and right minimal since $\pi$ is a
  projective cover.
\end{proof}

\section{Auslander varieties}

In \cite[\S6]{Ri2013}, Ringel pointed out the geometric nature of the
correspondence between submodules of $\Hom_\A(\C,Y)$ and right
$\C$-determined morphisms ending at an object $Y$ (see
Theorem~\ref{th:mor}).  More precisely, he defines under suitable
assumptions for a dimension vector $\mathbf d$ the \emph{Auslander
  variety} \[\Gr_{\mathbf d}(\Hom_\A(\C,Y))\] as an algebraic variety
given by all submodules of $\Hom_\A(\C,Y)$ with dimension vector
$\mathbf d$; it parametrizes right $\C$-determined morphisms ending
at $Y$. We illustrate this by giving an example.

Let $k$ be a field and fix a projective variety
$\bbX=V(f_1,\ldots,f_r)$ given by homogeneous polynomials $f_i\in
k[x_0,\ldots,x_n]$ of degree at most $p$.  This variety can be
realised as a \emph{quiver Grassmannian} \cite[\S2]{Hi1996}, and we
follow the exposition in \cite{Hu2014}.  Let $\La=\La_{n,p}$ denote
the \emph{Beilinson algebra } given by the path algebra of the
following quiver
\[\begin{tikzcd}[row sep=small]
\scriptstyle{p}
\arrow[yshift=1.0ex]{r}{x_0} \arrow[yshift=-1.0ex]{r}{\scriptstyle\cdots}[swap]{x_n}
&
{\scriptstyle p-1} \quad\cdots\quad {\scriptstyle 2}
\arrow[yshift=1.0ex]{r}{x_0}
\arrow[yshift=-1.0ex]{r}{\scriptstyle\cdots}[swap]{x_n}
&\scriptstyle{1}
\arrow[yshift=1.0ex]{r}{x_0} \arrow[yshift=-1.0ex]{r}{\scriptstyle\cdots}[swap]{x_n}
&{\scriptstyle 0}
\end{tikzcd}\] modulo all relations of the form $x_ix_j-x_jx_i$.  Each
homogeneous polynomial $f\in k[x_0,\ldots,x_n]$ of degree $d$ yields
$p-d+1$ elements of $\La$ represented by paths ending at vertices
$0,1,\ldots,p-d$, and we denote by $\langle f_1,\ldots, f_r\rangle$
the ideal of $\La$ generated by all occurences of each $f_i$. Consider
the indecomposable injective $\La/\langle f_1,\ldots, f_r\rangle$-
module $I(0)$ corresponding to the vertex $0$, but viewed as $\La$-module.
 Then $\bbX$ is isomorphic to the
variety of subrepresentations of $I(0)$ with
dimension vector $\mathbf d=(1,\ldots,1)$. It follows
that \[\bbX\cong\Gr_{\mathbf d}(\Hom_\A(\C,Y))\] for $\A= \Mod\La$,
$\C=\La$, and $Y=I(0)$.

\end{appendix}

\end{document}